\documentclass{amsproc}%
\usepackage{graphicx}
\usepackage{amscd}
\usepackage{amsmath}
\usepackage{amsfonts}
\usepackage{amssymb}%
\theoremstyle{plain}

\numberwithin{equation}{section}

\begin{document}
\title[A Cayley-Hamilton trace identity for $2\times2$ matrices]{A Cayley-Hamilton trace identity for $2\times2$ matrices over Lie-solvable rings}
\author{Johan Meyer}
\address{Department of Mathematics and Applied Mathematics, University of the Free
State, PO Box 339, Bloemfontein 9300, South Africa}
\email{MeyerJH.sci@ufs.ac.za}
\author{Jen\H{o} Szigeti}
\address{Institute of Mathematics, University of Miskolc, Miskolc, Hungary 3515}
\email{jeno.szigeti@uni-miskolc.hu}
\author{Leon van Wyk}
\address{Department of Mathematical Sciences, Stellenbosch University\\
P/Bag X1, Matieland 7602, Stellenbosch, South Africa }
\email{LvW@sun.ac.za}
\thanks{The first and third authors were supported by the National Research Foundation
of South Africa under Grant No.~UID 72375. Any opinion, findings and
conclusions or recommendations expressed in this material are those of the
authors and therefore the National Research Foundation does not accept any
liability in regard thereto.}
\thanks{As far as the second author is concerned, this research was carried out as
part of the TAMOP-4.2.1.B-10/2/KONV-2010-0001 project with support by the
European Union, co-financed by the European Social Fund.}
\thanks{The authors thank P. N. Anh and L. Marki for fruitful consultations.}

\begin{abstract}
We exhibit a Cayley-Hamilton trace identity for $2\times2$ matrices with
entries in a ring $R$ satisfying $[[x,y],[x,z]]=0$ and $\frac{1}{2}\in R$.

\end{abstract}
\maketitle

\noindent1. INTRODUCTION

\bigskip

The Cayley-Hamilton theorem and the corresponding trace identity play a
fundamental role in proving classical results about the polynomial and trace
identities of the $n\times n$ matrix algebra $M_{n}(K)$ over a field $K$ (see
[2] and [3]).\ In case of char$(K)=0$, Kemer's pioneering work (see [5]) on
the T-ideals of associative algebras revealed the importance of the identities
satisfied by the $n\times n$ matrices over the Grassmann (exterior) algebra%
\[
E=K\left\langle v_{1},v_{2},...,v_{r},...\mid v_{i}v_{j}+v_{j}v_{i}=0\text{
for all }1\leq i\leq j\right\rangle
\]
generated by the infinite sequence of anticommutative indeterminates
$(v_{i})_{i\geq1}$.

For $n\times n$ matrices over a Lie-nilpotent ring $R$ satisfying the
polynomial identity%
\[
\lbrack\lbrack\lbrack...[[x_{1},x_{2}],x_{3}],...],x_{m}],x_{m+1}]=0
\]
(with $[x,y]=xy-yx$), a Cayley-Hamilton identity of degree $n^{m}$ (with left-
or right-sided scalar coefficients) was found in [6]. Since $E$ is
Lie-nilpotent with $m=2$, the above mentioned Cayley-Hamilton identity for a
matrix $A\in M_{n}(E)$ is of degree$~n^{2}$.

In [1] Domokos presented a slightly modified version of this identity in which
the coefficients are invariant under the conjugate action of $GL_{n}(K)$. For
a matrix $A\in M_{2}(E)$ he obtained the trace identity%
\[
A^{4}\!-\!2\text{tr}(A)\!A^{3}\!+\!\left(  \!2\text{tr}^{2}(A)\!-\!\text{tr}%
(A^{2})\!\right)  \!A^{2}\!+\!\left(  \frac{1}{2}\!\text{tr}(A)\text{tr}%
(A^{2})\!+\!\frac{1}{2}\!\text{tr}(A^{2})\text{tr}(A)\!-\!\text{tr}%
^{3}(A)\!\right)  \!A\!+
\]%
\[
\frac{1}{4}\!\!\left(  \!\!\text{tr}^{4}(A)\!+\!\text{tr}^{2}(A^{2}%
)\!-\!\frac{5}{2}\!\text{tr}^{2}(A)\!\text{tr}(A^{2})\!+\!\frac{1}%
{2}\!\text{tr}(A^{2})\!\text{tr}^{2}(A)\!-2\text{tr}(A^{3})\!\text{tr}%
(A)\!+\!2\text{tr}(A)\!\text{tr}(A^{3})\!\!\right)  \!\!I\!\!=\!\!0,
\]
where $I$\ is the identity matrix and tr$(A)$ denotes the trace of $A$. A
similar identity with right coefficients also holds for $A$. Here $E$ can be
replaced by any ring $R$ which is Lie-nilpotent of index $2$.

The identity $[x,y][x,z]=0$ is a consequence of Lie-nilpotency of index $2$
(see [4]), as is obviously $[[x,y],[x,z]]=0$. The first aim of the present
paper is to provide an example of an algebra satisfying $[[x,y],[u,v]]=0$, but
neither $[x,y][u,v]=0$ nor $[[x,y],z]=0$. Since the above mentioned trace
identity cannot be used for matrices over such an algebra, our second purpose
is to exhibit a new trace identity of the same kind (of degree $4$ in $A$) for
a matrix $A$ in $M_{2}(R)$, where $R$ is any ring satisfying the identity%
\[
\lbrack\lbrack x,y],[x,z]]=0
\]
and $\frac{1}{2}\in R$. We note that a ring satisfying $[[x,y],[u,v]]=0$ is
called Lie-solvable of index $2$.

From now onward $R$ and $S$ are rings with $1$. In Section 2 we consider the
ring $U_{3}^{\ast}(R)$\ of upper triangular $3\times3$ matrices with equal
diagonal entries over $R$. First we observe that $U_{3}^{\ast}(R)$\ is never
commutative. We prove that if $R$ is commutative then the algebra $U_{3}%
^{\ast}(R)$\ satisfies the identities $[x,y][u,v]=0$ and $[[x,y],z]=0$.
However, for a non-commutative $R$ we show that the ring $U_{3}^{\ast}(R)$
never satisfies any of the identities $[x,y][u,v]=0$ and $[[x,y],z]=0$.

The main result in Section 2 states that if $S$ satisfies the identities
$[x,y][u,v]=0$ and $[[x,y],z]=0$, then the matrix ring $U_{3}^{\ast}(S)$\ is
Lie-solvable of index $2$. It follows that if $R$ is commutative, then
$U_{3}^{\ast}(U_{3}^{\ast}(R))$ is an example of an algebra satisfying
$[[x,y],[u,v]]=0$, but neither $[x,y][u,v]=0$ nor $[[x,y],z]=0$.

Section 3 is entirely devoted to the construction of our Cayley-Hamilton trace identity.

\bigskip

\noindent2. A PARTICULAR\ LIE-SOLVABLE MATRIX\ ALGEBRA

\bigskip

Since%
\[
E_{1,2},E_{2,3}\in U_{3}^{\ast}(R)=\left\{  \left[
\begin{array}
[c]{ccc}%
a & b & c\\
0 & a & d\\
0 & 0 & a
\end{array}
\right]  \mid a,b,c,d\in R\right\}
\]
and $E_{1,2}E_{2,3}=E_{1,3}\neq0=E_{2,3}E_{1,2}$, the ring $U_{3}^{\ast}(R)$
is never commutative. Any element of $U_{3}^{\ast}(R)$ can be written as
$aI+X$, where $X$ is strictly upper triangular. We note that $XYZ=0$ for
strictly upper triangular $3\times3$ matrices. If $R$ is commutative, then
$aI$ is central in $U_{3}^{\ast}(R)$ (of course, also in $M_{3}(R)$),
$[aI+X,bI+Y]=[X,Y]$ for all $a,b\in R$ and so $U_{3}^{\ast}(R)$ satisfies all
polynomial identities in which each summand is a product of certain (possibly
iterated) commutators. For example,%
\[
\lbrack x,y][u,v]=0\text{ and }[[x,y],z]=0
\]
are typical such identities for $U_{3}^{\ast}(R)$. If $R$ is non-commutative,
say $[r,s]\neq0$ for some $r,s\in R$, then for $x=rI$, $y=sE_{1,2}$,
$u=E_{2,2}$, $v=z=E_{2,3}$ in $U_{3}^{\ast}(R)$ we have%
\[
\lbrack x,y][u,v]=[[x,y],z]=[r,s]E_{1,3}\neq0.
\]

\bigskip

\noindent\textbf{2.1.Theorem.}\textit{ If }$S$\textit{ satisfies
}$[x,y][u,v]=0$\textit{ and }$[[x,y],z]=0$\textit{, then }$U_{3}^{\ast}%
(S)$\textit{ satisfies }$[[x,y],[u,v]]=0$\textit{.}

\bigskip

\noindent\textbf{Proof.} Using the matrices%
\[
x=\left[
\begin{array}
[c]{ccc}%
a & b & c\\
0 & a & d\\
0 & 0 & a
\end{array}
\right]  \text{ and }y=\left[
\begin{array}
[c]{ccc}%
e & f & g\\
0 & e & h\\
0 & 0 & e
\end{array}
\right]
\]
in $U_{3}^{\ast}(S)$, a straightforward calculation gives that%
\[
\lbrack x,y]=\left[
\begin{array}
[c]{ccc}%
\lbrack a,e] & [a,f]+[b,e] & [a,g]+[c,e]+(bh-fd)\\
0 & [a,e] & [a,h]+[d,e]\\
0 & 0 & [a,e]
\end{array}
\right]  =[a,e]I+C+\alpha E_{1,3},
\]
where $\alpha=bh-fd$ and $C$ is a strictly upper triangular matrix with
entries in $[S,S]$ (the additive subgroup of $S$ generated by all
commutators). Now $[[a,e],s]=0$ for all $s\in S$, hence $[a,e]I$ is central in
$U_{3}^{\ast}(S)$ (also in $M_{3}(S)$). Thus we have%
\[
\lbrack\lbrack x,y],[u,v]]=[[a,e]I+C+\alpha E_{1,3},[a^{\prime},e^{\prime
}]I+C^{\prime}+\alpha^{\prime}E_{1,3}]=[C+\alpha E_{1,3},C^{\prime}%
+\alpha^{\prime}E_{1,3}]=0
\]
because of $(C+\alpha E_{1,3})(C^{\prime}+\alpha^{\prime}E_{1,3})=(C^{\prime
}+\alpha^{\prime}E_{1,3})(C+\alpha E_{1,3})=0$. Indeed, $CC^{\prime}%
=C^{\prime}C=0$ is a consequence of $C,C^{\prime}\in M_{3}([S,S])$ and of
$[x,y][u,v]=0$ in $S$, and $CE_{1,3}=E_{1,3}C=C^{\prime}E_{1,3}=E_{1,3}%
C^{\prime}=0$ follows from the fact that $C$ and $C^{\prime}$ are strictly
upper triangular.$\square$

\bigskip

\noindent\textbf{2.2.Corollary.}\textit{ If }$R$\textit{ is commutative, then
the algebra }$U_{3}^{\ast}(U_{3}^{\ast}(R))$\textit{ satisfies}

\noindent$\lbrack\lbrack x,y],[u,v]]=0$\textit{, but neither }$[x,y][u,v]=0$%
\textit{ nor }$[[x,y],z]=0$\textit{.}

\bigskip

\noindent3. MATRICES\ WITH\ COMMUTATOR\ ENTRIES

\bigskip

The following can be considered as the \textquotedblleft
real\textquotedblright\ $2\times2$ Cayley-Hamilton trace identity.

\bigskip

\noindent\textbf{3.1.Proposition.}\textit{ If }$\frac{1}{2}\in R$\textit{ and
}$A=[a_{ij}]\in M_{2}(R)$\textit{, then}%
\[
A^{2}\!-\!\text{tr}(A)A\!+\!\frac{1}{2}\!(\!\text{tr}^{2}(A)\!-\!\text{tr}%
(A^{2})\!)I\!=\!\left[  \!%
\begin{array}
[c]{cc}%
\frac{1}{2}\![\!a_{11},\!a_{22}\!]\!+\!\frac{1}{2}\![\!a_{12},\!a_{21}\!] &
[a_{12},a_{22}]\\
\lbrack a_{21},a_{11}] & -\frac{1}{2}\![\!a_{11},\!a_{22}\!]\!-\!\frac{1}%
{2}\![\!a_{12},\!a_{21}\!]
\end{array}
\!\right]  .
\]
\noindent\textbf{Proof.} A straightforward computation suffices.$\square$

\bigskip

\noindent\textbf{3.2.Corollary.}\textit{ If }$\frac{1}{2}\in R$\textit{ and
}$B=[b_{ij}]\in M_{2}(R)$\textit{ with }tr$(B)=0$\textit{, then}%
\[
B^{2}-\frac{1}{2}\text{tr}(B^{2})I=\left[
\begin{array}
[c]{cc}%
\frac{1}{2}[b_{12},b_{21}] & -[b_{12},b_{11}]\\
\lbrack b_{21},b_{11}] & -\frac{1}{2}[b_{12},b_{21}]
\end{array}
\right]  .
\]
\noindent\textbf{Proof.} Since $b_{22}=-b_{11}$, we have $[b_{11},b_{22}]=0$
and $[b_{12},b_{22}]=-[b_{12},b_{11}]$. Thus the formula in Proposition 3.1
immediately gives the identity for $B$.$\square$

\bigskip

\noindent\textbf{3.3.Theorem.}\textit{ If }$\frac{1}{2}\in R$\textit{ and }%
$R$\textit{\ satisfies }$[[x,y],[x,z]]=0$\textit{, then}%
\[
\left(  C^{2}-\frac{1}{2}\text{tr}(C^{2})I\right)  ^{2}-\frac{1}{2}%
\text{tr}\left(  (C^{2}-\frac{1}{2}\text{tr}(C^{2})I)^{2}\right)  I=0
\]
\textit{for all }$C\in M_{2}(R)$\textit{ with }tr$(C)=0$\textit{.}

\bigskip

\noindent\textbf{Proof.} Take $C=[c_{ij}]$. In view of Corollary 3.2 we have%
\[
C^{2}-\frac{1}{2}\text{tr}(C^{2})I=\left[
\begin{array}
[c]{cc}%
\frac{1}{2}[c_{12},c_{21}] & -[c_{12},c_{11}]\\
\lbrack c_{21},c_{11}] & -\frac{1}{2}[c_{12},c_{21}]
\end{array}
\right]  .
\]
Since tr$(C^{2}-\frac{1}{2}$tr$(C^{2})I)=0$, the repeated application of
Corollary 3.2 to $B=C^{2}-\frac{1}{2}$tr$(C^{2})I$ gives that%
\[
\left(  C^{2}-\frac{1}{2}\text{tr}(C^{2})I\right)  ^{2}-\frac{1}{2}%
\text{tr}\left(  (C^{2}-\frac{1}{2}\text{tr}(C^{2})I)^{2}\right)  I=
\]%
\[
=\frac{1}{2}\left[
\begin{array}
[c]{cc}%
-[[c_{12},c_{11}],[c_{21},c_{11}]] & [[c_{12},c_{11}],[c_{12},c_{21}]]\\
\lbrack\lbrack c_{21},c_{11}],[c_{12},c_{21}]] & [[c_{12},c_{11}%
],[c_{21},c_{11}]]
\end{array}
\right]  .
\]
Now we have%
\[
\lbrack\lbrack c_{12},c_{11}],[c_{21},c_{11}]]=[[c_{11},c_{12}],[c_{11}%
,c_{21}]]
\]
and%
\[
\lbrack\lbrack c_{21},c_{11}],[c_{12},c_{21}]]=-[[c_{21},c_{11}],[c_{21}%
,c_{12}]].
\]
Thus each entry of the above $2\times2$ matrix is of the form $\pm
\lbrack\lbrack x,y],[x,z]]=0$ and the desired identity follows.$\square$

\bigskip

In Corollaries 3.4 - 3.6 we assume that $\frac{1}{2}\in R$ and $R$
satisfies\textit{ }$[[x,y],[x,z]]=0$.

\bigskip

\noindent\textbf{3.4.Corollary.}\textit{ If }$C\in M_{2}(R)$\textit{ with
}tr$(C)=0$\textit{,} \textit{then }%
\[
C^{4}-\frac{1}{2}\text{tr}(C^{2})C^{2}-\frac{1}{2}C^{2}\text{tr}(C^{2}%
)+\frac{1}{2}\left(  \text{tr}^{2}(C^{2})-\text{tr}(C^{4})\right)  I=0.
\]
\noindent\textbf{Proof.} Clearly,%
\[
\left(  C^{2}-\frac{1}{2}\text{tr}(C^{2})I\right)  ^{2}=C^{4}-\frac{1}%
{2}\text{tr}(C^{2})C^{2}-\frac{1}{2}C^{2}\text{tr}(C^{2})+\frac{1}{4}%
\text{tr}^{2}(C^{2})I
\]
and%
\[
\text{tr}\left(  (C^{2}\!-\!\frac{1}{2}\text{tr}(C^{2})I)^{2}\right)
\!=\!\text{tr}(C^{4})\!-\!\frac{1}{2}\text{tr}(\text{tr}(C^{2})C^{2}%
)\!-\!\frac{1}{2}\text{tr}(C^{2}\text{tr}(C^{2}))\!+\!\frac{1}{4}%
\text{tr}(\text{tr}^{2}(C^{2})I)\!=
\]%
\[
=\text{tr}(C^{4})-\frac{1}{2}\text{tr}^{2}(C^{2})-\frac{1}{2}\text{tr}%
^{2}(C^{2})+\frac{1}{2}\text{tr}^{2}(C^{2})=\text{tr}(C^{4})-\frac{1}%
{2}\text{tr}^{2}(C^{2}).
\]
Thus we have%
\[
\left(  C^{2}-\frac{1}{2}\text{tr}(C^{2})I\right)  ^{2}-\frac{1}{2}%
\text{tr}\left(  (C^{2}-\frac{1}{2}\text{tr}(C^{2})I)^{2}\right)  I=
\]%
\[
=C^{4}-\frac{1}{2}\text{tr}(C^{2})C^{2}-\frac{1}{2}C^{2}\text{tr}(C^{2}%
)+\frac{1}{4}\text{tr}^{2}(C^{2})I-\frac{1}{2}\left(  \text{tr}(C^{4}%
)-\frac{1}{2}\text{tr}^{2}(C^{2})\right)  I=
\]%
\[
=C^{4}-\frac{1}{2}\text{tr}(C^{2})C^{2}-\frac{1}{2}C^{2}\text{tr}(C^{2}%
)+\frac{1}{2}\left(  \text{tr}^{2}(C^{2})-\text{tr}(C^{4})\right)  I.\square
\]

\bigskip

\noindent\textbf{3.5.Corollary.}\textit{ If }$C\in M_{2}(R)$\textit{ with
}tr$(C)=$ tr$(C^{2})=$ tr$(C^{4})=0$\textit{, then }$C^{4}=0$\textit{.}

\bigskip

\noindent\textbf{3.6.Corollary.} \textit{If }$A\in M_{2}(R)$\textit{ is
arbitrary}, \textit{then }%
\[
\left(  A-\frac{1}{2}\text{tr}(A)I\right)  ^{4}-\frac{1}{2}\text{tr}\left(
(A-\frac{1}{2}\text{tr}(A)I)^{2}\right)  \left(  A-\frac{1}{2}\text{tr}%
(A)I\right)  ^{2}-
\]%
\[
\frac{1}{2}\left(  A-\frac{1}{2}\text{tr}(A)I\right)  ^{2}\text{tr}\left(
(A-\frac{1}{2}\text{tr}(A)I)^{2}\right)  +
\]%
\[
\frac{1}{2}\left(  \text{tr}^{2}\left(  (A-\frac{1}{2}\text{tr}(A)I)^{2}%
\right)  -\text{tr}\left(  (A-\frac{1}{2}\text{tr}(A)I)^{4}\right)  \right)
I=0.
\]

\bigskip

\noindent\textbf{Proof.} Take $C=A-\frac{1}{2}$tr$(A)I$, then tr$(C)=$
tr$(A-\frac{1}{2}$tr$(A)I)=0$ and the application of Corollary 3.4 gives the
identity.$\square$

\bigskip

\noindent\textbf{3.7.Theorem.} \textit{If }$\frac{1}{2}\in R$\textit{ and }%
$R$\textit{\ is a Lie-solvable ring satisfying }$[[x,y],[x,z]]=0$\textit{,
then for all }$A\in M_{2}(R)$\textit{ we have}%
\[
A^{4}-\frac{1}{2}A^{2}\text{tr}(A)A-\frac{1}{2}A\text{tr}(A)A^{2}-\frac{1}%
{2}A^{3}\text{tr}(A)-\frac{1}{2}\text{tr}(A)A^{3}+\frac{1}{2}A^{2}%
\text{tr}^{2}(A)+\frac{1}{2}\text{tr}^{2}(A)A^{2}-
\]%
\[
\frac{1}{2}A^{2}\text{tr(}A^{2})-\frac{1}{2}\text{tr(}A^{2})A^{2}+\frac{1}%
{4}A\text{tr}(A)A\text{tr}(A)+\frac{1}{4}\text{tr}(A)A\text{tr}(A)A+
\]%
\[
\frac{1}{4}\text{tr}(A)A^{2}\text{tr}(A)+\frac{1}{4}A\text{tr}^{2}%
(A)A-\frac{1}{4}\text{tr}(A)A\text{tr}^{2}(A)-\frac{1}{4}\text{tr}%
^{2}(A)A\text{tr}(A)+
\]%
\[
\frac{1}{4}\text{tr}(A)A\text{tr(}A^{2})+\frac{1}{4}\text{tr(}A^{2}%
)A\text{tr}(A)-\frac{1}{4}A\text{tr}^{3}(A)-\frac{1}{4}\text{tr}^{3}(A)A+
\]%
\[
\frac{1}{4}A\text{tr}(A)\text{tr(}A^{2})+\frac{1}{4}\text{tr(}A^{2}%
)\text{tr}(A)A-\frac{1}{2}\text{tr}^{2}(A)\text{tr(}A^{2})I-\frac{1}%
{2}\text{tr(}A^{2})\text{tr}^{2}(A)I+
\]%
\[
\frac{1}{2}\text{tr}^{2}\text{(}A^{2})I\!+\!\frac{1}{4}\text{tr}\left(
A^{2}\text{tr}(A)A\right)  I\!+\!\frac{1}{4}\text{tr}(A\text{tr}%
(A)A^{2})I\!+\!\frac{1}{4}\text{tr}(A^{3})\text{tr}(A)I\!+\!\frac{1}%
{4}\text{tr}(A)\text{tr}(A^{3})I\!-
\]%
\[
\frac{1}{8}\text{tr}(A)\text{tr}(A\text{tr}(A)A)I\!\!-\!\!\frac{1}{8}%
\text{tr}(A\text{tr}(A)A)\text{tr}(A)I\!\!-\!\!\frac{1}{8}\text{tr}%
(A\text{tr}^{2}(A)A)I\!\!-\!\!\frac{1}{8}\text{tr}(A)\text{tr(}A^{2}%
)\text{tr}(A)I\!+
\]%
\[
\frac{1}{2}\text{tr}^{4}(A)I-\frac{1}{2}\text{tr(}A^{4})I=0\text{.}%
\]

\bigskip

\noindent\textbf{Proof.} Take $C=A-\frac{1}{2}$tr$(A)I$, then tr$(C)=$%
tr$(A-\frac{1}{2}$tr$(A)I)=0$. We have%
\[
C^{2}-\frac{1}{2}\text{tr}(C^{2})I=(A-\frac{1}{2}\text{tr}(A)I)^{2}-\frac
{1}{2}\text{tr}\left(  (A-\frac{1}{2}\text{tr}(A)I)^{2}\right)  I=
\]%
\[
=\!A^{2}\!-\!\frac{1}{2}\text{tr}(A)A\!-\!\frac{1}{2}A\text{tr}(A)\!+\!\frac
{1}{4}\text{tr}^{2}(A)I\!-\!\frac{1}{2}\text{tr}\!\!\left(  \!\!A^{2}%
\!-\!\frac{1}{2}\text{tr}(A)A\!-\!\frac{1}{2}A\text{tr}(A)\!+\!\frac{1}%
{4}\text{tr}^{2}(A)I\!\!\right)  \!\!I\!=
\]%
\[
=\!A^{2}\!-\!\frac{1}{2}\text{tr}(A)A\!-\!\frac{1}{2}A\text{tr}(A)\!+\!\frac
{1}{4}\text{tr}^{2}(A)I\!-\!\frac{1}{2}\left(  \!\!\text{tr(}A^{2}%
)\!-\!\frac{1}{2}\text{tr}^{2}(A)\!-\!\frac{1}{2}\text{tr}^{2}(A)\!+\!\frac
{1}{2}\text{tr}^{2}(A)\!\!\right)  \!\!I\!=
\]%
\[
=A^{2}-\frac{1}{2}\text{tr}(A)A-\frac{1}{2}A\text{tr}(A)+\frac{1}{2}%
\text{tr}^{2}(A)I-\frac{1}{2}\text{tr(}A^{2})I
\]
and
\[
\left(  C^{2}-\frac{1}{2}\text{tr}(C^{2})I\right)  ^{2}=\left(  A^{2}-\frac
{1}{2}\text{tr}(A)A-\frac{1}{2}A\text{tr}(A)+\frac{1}{2}\text{tr}%
^{2}(A)I-\frac{1}{2}\text{tr(}A^{2})I\right)  ^{2}=
\]%
\[
=A^{4}-\frac{1}{2}A^{2}\text{tr}(A)A-\frac{1}{2}A^{3}\text{tr}(A)+\frac{1}%
{2}A^{2}\text{tr}^{2}(A)-\frac{1}{2}A^{2}\text{tr(}A^{2})-\frac{1}{2}%
\text{tr}(A)A^{3}+
\]%
\[
\frac{1}{4}\text{tr}(A)A\text{tr}(A)A+\frac{1}{4}\text{tr}(A)A^{2}%
\text{tr}(A)-\frac{1}{4}\text{tr}(A)A\text{tr}^{2}(A)+\frac{1}{4}%
\text{tr}(A)A\text{tr(}A^{2})-
\]%
\[
\frac{1}{2}A\text{tr}(A)A^{2}+\frac{1}{4}A\text{tr}^{2}(A)A+\frac{1}%
{4}A\text{tr}(A)A\text{tr}(A)-\frac{1}{4}A\text{tr}^{3}(A)+\frac{1}%
{4}A\text{tr}(A)\text{tr(}A^{2})+
\]%
\[
\frac{1}{2}\text{tr}^{2}(A)A^{2}-\frac{1}{4}\text{tr}^{3}(A)A-\frac{1}%
{4}\text{tr}^{2}(A)A\text{tr}(A)+\frac{1}{4}\text{tr}^{4}(A)I-\frac{1}%
{4}\text{tr}^{2}(A)\text{tr(}A^{2})I-
\]%
\[
\frac{1}{2}\text{tr(}A^{2})A^{2}+\frac{1}{4}\text{tr(}A^{2})\text{tr}%
(A)A+\frac{1}{4}\text{tr(}A^{2})A\text{tr}(A)-\frac{1}{4}\text{tr(}%
A^{2})\text{tr}^{2}(A)I+\frac{1}{4}\text{tr}^{2}\text{(}A^{2})I.
\]
Thus we obtain that%
\[
\text{tr}\left(  C^{2}-\frac{1}{2}\text{tr}(C^{2})I\right)  ^{2}%
=\text{tr}(A^{4})-\frac{1}{2}\text{tr}\left(  A^{2}\text{tr}(A)A\right)
-\frac{1}{2}\text{tr}(A^{3})\text{tr}(A)+
\]%
\[
\frac{1}{2}\text{tr}(A^{2})\text{tr}^{2}(A)-\frac{1}{2}\text{tr}%
(A^{2})\text{tr(}A^{2})-\frac{1}{2}\text{tr}(A)\text{tr}(A^{3})+
\]%
\[
\frac{1}{4}\text{tr}(A)\text{tr}(A\text{tr}(A)A)+\frac{1}{4}\text{tr}%
(A)\text{tr(}A^{2})\text{tr}(A)-
\]%
\[
\frac{1}{4}\text{tr}(A)\text{tr}(A)\text{tr}^{2}(A)+\frac{1}{4}\text{tr}%
(A)\text{tr}(A)\text{tr(}A^{2})-\frac{1}{2}\text{tr}(A\text{tr}(A)A^{2}%
)+\frac{1}{4}\text{tr}(A\text{tr}^{2}(A)A)+
\]%
\[
\frac{1}{4}\text{tr}(A\text{tr}(A)A)\text{tr}(A)-\frac{1}{4}\text{tr}%
(A)\text{tr}^{3}(A)+\frac{1}{4}\text{tr}(A)\text{tr}(A)\text{tr(}A^{2}%
)+\frac{1}{2}\text{tr}^{2}(A)\text{tr}(A^{2})-
\]%
\[
\frac{1}{4}\text{tr}^{3}(A)\text{tr}(A)-\frac{1}{4}\text{tr}^{2}%
(A)\text{tr}(A)\text{tr}(A)+
\]%
\[
\frac{1}{2}\text{tr}^{4}(A)-\frac{1}{2}\text{tr}^{2}(A)\text{tr(}A^{2}%
)-\frac{1}{2}\text{tr(}A^{2})\text{tr}(A^{2})+
\]%
\[
\frac{1}{4}\text{tr(}A^{2})\text{tr}(A)\text{tr}(A)+\frac{1}{4}\text{tr(}%
A^{2})\text{tr}(A)\text{tr}(A)-\frac{1}{2}\text{tr(}A^{2})\text{tr}%
^{2}(A)+\frac{1}{2}\text{tr}^{2}\text{(}A^{2})=
\]

\[
=\text{tr}(A^{4})-\frac{1}{2}\text{tr}\left(  A^{2}\text{tr}(A)A\right)
-\frac{1}{2}\text{tr}(A\text{tr}(A)A^{2})-
\]%
\[
\frac{1}{2}\text{tr}(A^{3})\text{tr}(A)-\frac{1}{2}\text{tr}(A)\text{tr}%
(A^{3})+\frac{1}{4}\text{tr}(A)\text{tr}(A\text{tr}(A)A)+\frac{1}{4}%
\text{tr}(A\text{tr}(A)A)\text{tr}(A)+
\]%
\[
\frac{1}{4}\text{tr}(A\text{tr}^{2}(A)A)+\frac{1}{4}\text{tr}(A)\text{tr(}%
A^{2})\text{tr}(A)+\frac{1}{2}\text{tr}^{2}(A)\text{tr(}A^{2})+\frac{1}%
{2}\text{tr(}A^{2})\text{tr}^{2}(A)-
\]%
\[
\frac{1}{2}\text{tr}^{2}\text{(}A^{2})-\frac{1}{2}\text{tr}^{4}(A).
\]
Now the calculation of%
\[
\left(  C^{2}-\frac{1}{2}\text{tr}(C^{2})I\right)  ^{2}-\frac{1}{2}%
\text{tr}\left(  (C^{2}-\frac{1}{2}\text{tr}(C^{2})I)^{2}\right)  I
\]
and the application of Theorem 3.3 yield the identity.$\square$

\bigskip

\noindent\textbf{Question.} Throughout this section we have used the identity
$[[x,y],[x,z]]=0$. We do not know whether this identity implies the
\textquotedblleft seemingly\textquotedblright\ stronger identity
$[[x,y],[u,v]]=0$ which plays an important role in Section 2.

\bigskip

Starting with a matrix $C\in M_{2}(R)$ such that tr$(C)=0$, define the
sequence $(C_{k})_{k\geq0}$ by the following recursion: $C_{0}=C$ and
\[
C_{k+1}=C_{k}^{2}-\frac{1}{2}\text{tr}(C_{k}^{2})I.
\]
Clearly, tr$(C_{k})=0$ for all $k\geq0$ and $C_{k}$ is a trace polynomial
expression of $C$. In view of Corollary 3.2, the entries of $C_{1}$ are of the
form $[x_{1},x_{2}]$. The repeated application of Corollary 3.2 (as it can be
seen in the proof of Theorem 3.3) and a straightforward induction show that
the (four) entries of $C_{k}$ are all of the form $[x_{1},x_{2},...,x_{2^{k}%
}]_{\text{solv}}$, where $[x_{1},x_{2}]_{\text{solv}}=[x_{1},x_{2}]$ and for
$i\geq1$ we take the Lie brackets as%
\[
\lbrack x_{1},x_{2},...,x_{2^{i+1}}]_{\text{solv}}=[[x_{1},x_{2},...,x_{2^{i}%
}]_{\text{solv}},[x_{2^{i}+1},x_{2^{i}+2},...,x_{2^{i}+2^{i}}]_{\text{solv}}].
\]
If $R$ satisfies the general identity%
\[
\lbrack x_{1},x_{2},...,x_{2^{k}}]_{\text{solv}}=0
\]
of Lie solvability, then $C_{k}=0$, whence we can derive a trace identity for
$C$. Thus the substitution $C=A-\frac{1}{2}$tr$(A)I$ gives a trace identity
for an arbitrary $A\in M_{2}(R)$.

\bigskip

\bigskip

\noindent REFERENCES

\bigskip

\begin{enumerate}
\item M. Domokos, \textit{Cayley-Hamilton theorem for }$2\times2$%
\textit{\ matrices over the Grassmann algebra,} J. Pure Appl. Algebra 133
(1998), 69-81.

\item V. Drensky, \textit{Free Algebras and PI-Algebras,} Springer-Verlag, 2000.

\item V. Drensky and E. Formanek, \textit{Polynomial Identity Rings,}
Birkh\"{a}user-Verlag, 2004.

\item S. A. Jennings, \textit{On rings whose associated Lie rings are
nilpotent,} Bull. Amer. Math. Soc. 53 (1947), 593-597.

\item A. R. Kemer,\textit{\ Ideals of Identities of Associative Algebras,}
Translations of Math. Monographs, Vol. 87 (1991), AMS, Providence, Rhode Island.

\item J. Szigeti,\textit{\ New determinants and the Cayley-Hamilton theorem
for matrices over Lie nilpotent rings,} Proc. Amer. Math. Soc. 125(8) (1997), 2245-2254.
\end{enumerate}

\end{document}